\documentclass[11pt]{article}
\usepackage{amsmath}
\usepackage{amsfonts}
\usepackage{bm}
\usepackage{graphicx}
\usepackage{esdiff}
\usepackage{url}

\newcommand{\be}{\begin{equation}}
\newcommand{\ee}{\end{equation}}
\newcommand{\ba}{\begin{array}}
\newcommand{\ea}{\end{array}}
\newcommand{\baa}{\left[\begin{array}}
\newcommand{\eaa}{\end{array}\right]}
\newcommand{\ds}{\displaystyle}
\newcommand{\beqa}{\begin{eqnarray}}
\newcommand{\eeqa}{\end{eqnarray}}
\newcommand{\bt}{\begin{tabular}}
\newcommand{\et}{\end{tabular}}

\newcommand{\bc}{\begin{center}}
\newcommand{\ec}{\end{center}}

\newtheorem{remark}{Remark}
\newtheorem{example}[remark]{Example}

\def\QED{\hfill \mbox{\rule[0pt]{1.5ex}{1.5ex}}}
\newcommand{\degpi}{\kappa}

\title{Semi-infinite optimization with sums of exponentials
via polynomial approximation}
\author{Bogdan Dumitrescu$^{\rm a,b}$,
 Bogdan C. \c{S}icleru$^{\rm a}$, Florin Avram$^{\rm c}$}

\begin{document}
\maketitle

\vspace*{20mm}
\begin{abstract}
We propose a general method for optimization with semi-infinite constraints
that involve a linear combination of functions, focusing on the
case of the exponential function.
Each function is lower and upper bounded on sub-intervals
by low-degree polynomials.
Thus, the constraints can be approximated with polynomial inequalities
that can be implemented with linear matrix inequalities.
Convexity is preserved, but the problem has now a finite number of
constraints.
We show how to take advantage of the properties of the exponential
function in order to build quickly accurate approximations.
The problem used for illustration is the least-squares fitting of
a positive sum of exponentials to an empirical density.
When the exponents are given, the problem is convex, but we also
give a procedure for optimizing the exponents.
Several examples show that the method is flexible, accurate and gives
better results than other methods for the investigated problems.
\end{abstract}

\vspace*{50mm}
\noindent $^{\rm a}$Department of Automatic Control and Computers,
University Politehnica of Bucharest, Spl.\ Independen\c{t}ei 313, Bucharest 060042, Romania.
E-mails: bogdan.dumitrescu, bogdan.sicleru@acse.pub.ro \\
$^{\rm b}$ Department of Signal Processing, Tampere University of Technology, Finland.
E-mail: bogdan.dumitrescu@tut.fi \\
$^{\rm c}$ Universit\'e de Pau et des Pays de l'Adour, France



\pagebreak
\section{Introduction}

\subsection{The problem}

The purpose of this paper is twofold.
Our main aim is to give a detailed procedure for solving
a class of semi-infinite programming (SIP) problems involving
the function
\be
f(t) = \sum_{i=1}^n p_i(t) e^{-\lambda_i t},
\label{sope}
\ee
where $0 \le \lambda_1 < \ldots < \lambda_n$
and $p_i(t) = \sum_{j=0}^{\degpi_i} \alpha_{i,j} t^j$ are polynomials,
having usually a low degree $\degpi_i$.
For the sake of abbreviation, we name SOPE (sum of polynomials times exponentials)
a function like (\ref{sope}).

In the same time, but with much less detail, we point out how the
proposed approach can be extended to a much more general
category of functions, with the exponentials from (\ref{sope})
replaced by arbitrary (typically elementary) functions. 

The main source of problems involving (\ref{sope})
is the modeling of random variables.
If $p_i(t)$ are all constants ($\degpi_i=0$), then $f(t)$ is
the probability density function (pdf) of a hyperexponential distribution.
If $n=1$ and $\alpha_{i,j}$ are nonnegative, then $f(t)$ is the pdf of an Erlang mixture.
If $n>1$ and $\alpha_{i,j}$ are nonnegative, then $f(t)$ is called
a hyper-Erlang density.

There are three types of optimization problems involving (\ref{sope})
that we tackle.
They all involve SIP constraints having the form
\be
f(t) \ge 0, \ \forall t \in [t_0, t_f],
\label{fpos}
\ee
where typical choices are $t_0 = 0$, $t_f = \infty$;
however, it is enough to take a sufficiently large $t_f$ instead
of infinity, due to the decay properties of the exponential.
For simplicity of reference, we give ad hoc names to these problems
and list them in increasing order of difficulty.

{\em Positivity check (SOPE-P).}
The simplest problem is: given a function (\ref{sope}),
decide if the inequality (\ref{fpos}) holds.
Such decision is necessary e.g.\ when $f(t)$ should be a pdf.

{\em Convex constraint (SOPE-C).}
More interesting are the cases where the parameters of the function
are variables of the optimization problem.
In particular, when the polynomial coefficients $\alpha_{i,j}$ are variable,
but the exponents $\lambda_i$ are given, the constraint (\ref{fpos})
is convex.
This type of problem is our main focus.
We will always assume that the number $n$ of exponentials and the degrees $\degpi_i$
of the polynomials are known; they can be selected using Information Theoretic
Criteria, but this is beyond the scope of this paper.

A typical problem is fitting a pdf to empirical data.
A random process is governed by an unknown pdf $h(t)$; from empirical observations,
we know the values $h_m = h(t_m) \ge 0$, for some times $t_m$, $m=1:M$, usually equidistant.
Assuming that (\ref{sope}) is an appropriate model, we want to fit it to the data
by solving the least squares problem
\be
\ba{ccl} \min & & \ds J(f) = \frac{1}{M} \sum_{m=1}^M w_m [ f(t_m) - h_m ]^2 \\
\mbox{s.t.} & & f(t) \ge 0, \ \forall t \in [t_0, t_f]
\ea
\label{Jsope}
\ee
where $w_m > 0$ are weights (we take implicitly $w_m=1$).
If the exponents $\lambda_i$ are given, the problem (\ref{Jsope}) belongs to convex SIP;
its difficulty comes from the infinite number of constraints hidden by (\ref{fpos}).
There is no simple way to express it equivalently using a finite number
of constraints like, for example, for polynomials.

With small modifications, the problem (\ref{Jsope}) can be posed not for pdfs,
but for cumulative density functions (cdf) or complementary cdfs, see later (\ref{JsopeF}).

{\em General constraint (SOPE-G).}
The most difficult problem is when the exponents $\lambda_i$ are also variable.
In this case an optimization problem like ({\ref{Jsope}}) is no longer convex. 
Although we will present a solution to this problem also,
it will have no guarantee of optimality.

An interesting particular case of the above problems is when $\degpi_i=0$, hence
the polynomials are reduced to constants and so
\be
f(t) = \sum_{i=1}^n \alpha_i e^{-\lambda_i t}.
\label{soe}
\ee
The name of the problems will be changed accordingly, by replacing SOPE with
SOE (sum of exponentials).

\subsection{Contribution and contents}

The approach we propose is based on a polynomial approximation
of the exponential function that allows the approximation
of a SOPE-C problem like (\ref{Jsope}) with a polynomial optimization
problem that can be expressed with linear matrix inequalities
and hence solved efficiently with semidefinite programming (SDP) methods
in friendly media like CVX \cite{CVX} and POS3POLY \cite{SiDu13oe}
or GloptiPoly \cite{HeLL09_GP3}.

We split the inequality (\ref{fpos}) in $K$ sub-intervals
\be
[t_0,t_f] = \bigcup_{k=1}^K [t_{k-1}, t_k]
\label{interv}
\ee
and impose positivity on each sub-interval.
We compute polynomials $\hat b_{ik}(t)$, $\check b_{ik}(t)$
such that
\be
\hat b_{ik}(t) \ge e^{-\lambda_i t} \ge \check b_{ik}(t),
\ \ \forall t \in [t_{k-1}, t_k].
\label{approxint}
\ee
Using such lower and upper approximations of the exponentials,
we impose sufficient positivity conditions on (\ref{sope}),
simultaneously on all sub-intervals.
The resulting approximation of SOPE-C will be detailed in Section \ref{sec:appol}.

One-sided approximations like in (\ref{approxint}) can be computed
using the results from \cite{BDV66} (see also \cite{DV68,Lew70}),
as reviewed later in Section \ref{sec:expapp}.
Since the derivatives of the exponential have constant sign,
a satisfactory approximation can be computed by solving a linear system.

It is intuitive that by playing with the degrees of the approximation
polynomials and the lengths of the sub-intervals, the approximation can
be in principle as good as desired, at the price of increased complexity.
This is true not only for the exponential, but for all smooth functions
on a finite interval.
Numerical accuracy is also an important issue, especially since the
approximating polynomials from (\ref{approxint}) are defined by
their coefficients.
However, since we use the inequalities (\ref{approxint}) to build
a lower approximation of $f(t)$ and thus positivity is not violated
(like when the exponentials would be approximated with a single
polynomial), the quality of the approximation should be 
comparable to the accuracy of the
SDP algorithm used for solving the transformed problem.

In a practical implementation, we have reached the conclusion that it is
enough to use relatively low degrees of the approximation polynomials,
for example 8 or 10, in order to reach reasonable accuracy and complexity.

Using this approach, we propose in Section \ref{sec:fit} a complete
iterative procedure for solving the SOE-G and SOPE-G problems (\ref{Jsope}).
In particular, we select initial values for the exponents $\lambda_i$
by searching sparse functions (\ref{sope}) whose exponents belong to 
an arithmetic progression.
For given exponents, the coefficients $\alpha_{i,j}$ are optimized using
the polynomial approximation suggested above.
Possibly better values of the exponents are sought through a
descent search, the coefficients are re-optimized, etc.

Although the procedure is based on simple ideas, the results
given in Section \ref{sec:res} show that it is competitive, giving
good results on several fitting problems studied previously.

\subsection{Previous work}

The idea of using polynomial approximations in optimization
is by no means new.
However, as far as we know, the particular combination of ideas
that lies at the foundation of our approach was not proposed.
We present below a few relevant works and show the distinctive features
of our method.

Most of the methods use a single polynomial for approximation on
the whole interval or on sub-intervals.
A leading example is the library Chebfun \cite{chebfunv4}, which allows
numerically-performed "symbolic" computation by actually replacing
{\em given} functions with polynomials that approximate them to
machine precision on the interval of interest.
However, due to numerical considerations, the polynomials are not defined
by their coefficients but by their values in Chebyshev nodes
and barycentric Lagrange interpolation is employed for
computing the function values.
Hence, Chebfun can solve a problem like SOPE-P by computing the
minimum of the function, but is not able to take advantage of convexity.

In \cite{DGN10}, functional optimization is performed by approximating the
unknown function $f(t)$ with a polynomial with unknown coefficients.
It can be seen as a {\em non-parametric} form of our problem,
where the values of $f(t)$ are sought for all $t$ in some interval.
Here, we are interested in the parameters $\alpha_{i,j}$ of a
function with known structure, hence our method is {\em parametric}.

A precursor of the above approach was presented in \cite{KEH06},
where the minimum of a given function was obtained
by computing instead the minimum of an interpolating polynomial.
Polynomial positivity, used for the minimum computation, is enforced
through SDP via interpolation \cite{LofP04}.

In somewhat the same vein of non-parametric methods, but with pdf
estimation as specific target,
functional approximation with polynomials of a pdf fitting given moments
is proposed in \cite{HLM12}.
A simpler idea appeared in \cite{Short05}, proposing piecewise polynomial
approximations for pdfs, with no interest on optimization and no positivity enforcement.
More refined approaches appeal to splines \cite{JAOT07,MTV08,AlPa13}
(the latter work summarizes previous work
of its authors, including the use of multidimensional splines);
positivity is enforced with an ad hoc method in the first work,
while in the others it is imposed on each spline section via SDP.

One-sided polynomial approximations are less current;
for example, in \cite{ShWa01}, they were used for relaxation
in a branch and bound process for global optimization.
The approximation regarded the optimization variables, not
independent functions, like here.
Polynomial optimization was solved via linear programming.

There is also an entire body of literature dedicated to the
estimation of the parameters of a pdf, with diverse applications.
Optimization is used more or less explicitly, but usually
in a relatively standard way.
We cite only a few methods connected to our work.
SDP was used in \cite{FHT06} for estimating pdfs that
are the product between a polynomial and a kernel with
few parameters, like the Erlang mixture; polynomial positivity
is imposed via SDP.
In \cite{Duf07}, the SOE-G problem is solved within the class
of exponents in arithmetic progression; however, the opportunity
is missed for transforming it into a polynomial problem, like
shown in Section \ref{sec:fit}, and hence use SDP; positivity
is imposed by rather ad hoc means.
In \cite{SOH12}, matrix exponential distributions are estimated;
after finding the exponents with a method for linear systems identification,
the coefficients are found by optimization, positivity being ensured
by a Budan-Fourier technique. 
Among other families of techniques, we mention expectation maximization
\cite{ANO96}.

Although this paper gives an incomplete image, our aim is visible.
When implemented not only for the exponential, but for more general
categories of functions, our method could be seen as a possible meeting
point between CVX and Chebfun, by extending CVX to SIP,
using polynomial approximations like Chebfun
and having a simple modus operandi like both of them.

\section{Transformation to a polynomial problem}
\label{sec:appol}

We discuss here the approximation of SOPE-C with an optimization
problem with polynomials, thus preserving convexity, but transforming
the semi-infinite constraint into a finite one.

Let us assume that the splitting (\ref{interv}) and the approximations
(\ref{approxint}) of the exponentials are available.
We show how to impose the positivity constraint (\ref{fpos})
on a generic sub-interval $[t_{k-1}, t_k]$, keeping in mind
that the conditions are imposed simultaneously on all sub-intervals.

\subsection{Known signs}

Assume first that the polynomials $p_i(t)$ have constant and known
signs on $[t_{k-1}, t_k]$.
This is usually possible only when solving SOPE-P, where $p_i(t)$ are fixed.
Let $I_+$ be the set of indices $i$ for which $p_i(t) \ge 0$ on $[t_{k-1}, t_k]$,
and $I_-$ be defined similarly for the negative case.
In this situation, we can write
\be
f(t) \ge \sum_{i \in I_+} p_i(t) \check b_{ik}(t)
+ \sum_{i \in I_-} p_i(t) \hat b_{ik}(t) = P_k(t).
\label{Psign}
\ee
Since $P_k(t)$ is a polynomial whose coefficients depend linearly on
those of $p_i(t)$, imposing $P_k(t) \ge 0$ on $[t_{k-1}, t_k]$ is easy
and can be inserted into any convex optimization problem involving $f(t)$.
How to transform a polynomial positivity constraint into a linear matrix
inequality is discussed in \cite{Nest00}.
However, such knowledge is not necessary when using
a library facilitating the manipulation of positive polynomials, like POS3POLY.

\subsection{Unknown signs}

In general, the coefficients and the signs of $p_i(t)$ are not known,
since the coefficients of these polynomials are variables
in the optimization problem SOPE-C.
In this case, it is impossible to build the approximation (\ref{Psign}).
However, we can replace it with
\be
P_k(t) = \sum_{i=1}^n [p_i(t) - \gamma_{ik}(t)] \hat b_{ik}(t)
+ \gamma_{ik}(t) \check b_{ik}(t),
\label{Pnosign}
\ee
with the extra conditions
\be
\ba{l}
p_i(t) - \gamma_{ik}(t) \le 0, \ \ \forall t \in [t_{k-1}, t_k] \\
\gamma_{ik}(t) \ge 0,
\ea
\label{picond}
\ee
where $\gamma_{ik}(t)$ is a polynomial, typically of the
same degree as $p_i(t)$.
If $\deg \gamma_{ik} < \deg p_i$, then we can consider the
alternative bounding polynomial
\be
P_k(t) = \sum_{i=1}^n \gamma_{ik}(t) \hat b_{ik}(t)
+ [p_i(t) - \gamma_{ik}(t)] \check b_{ik}(t),
\label{Pnosign1}
\ee
with the extra conditions
\be
\ba{l}
p_i(t) - \gamma_{ik}(t) \ge 0, \ \ \forall t \in [t_{k-1}, t_k] \\
\gamma_{ik}(t) \le 0.
\ea
\label{picond1}
\ee

If $\deg \gamma_{ik} \ge \deg p_i$, then (\ref{Pnosign})-(\ref{picond})
and (\ref{Pnosign1})-(\ref{picond1}) are equivalent;
otherwise, they are different, but it is hard to give general rules
for choosing one over the other.
In both cases, it is clear that $f(t) \ge P_k(t)$ and so
$P_k(t) \ge 0$ is a sufficient condition for $f(t) \ge 0$.
From now on we will work only with (\ref{Pnosign})-(\ref{picond}).

So, the SOPE-C problem (\ref{Jsope}) is approximated with
\be
\ba{ccl} \min & & \ds J(f) \\
\mbox{s.t.} & & \sum_{i=1}^n [p_i(t) - \gamma_{ik}(t)] \hat b_{ik}(t)
                + \gamma_{ik}(t) \check b_{ik}(t) \ge 0,
				\ \ \forall t \in [t_{k-1}, t_k], \ k=1:K \\
            & & \left.
                \ba{l} p_i(t) - \gamma_{ik}(t) \le 0 \\
                  \gamma_{ik}(t) \ge 0
                \ea \right\} \ \ \forall t \in [t_{k-1}, t_k], \ k=1:K, \ i=1:n
\ea
\label{SOPE-C_app}
\ee
Note that the criterion is unchanged.
The price for generality is the apparition of the new variable polynomials
$\gamma_{ik}(t)$ and of the new constraints (\ref{picond}).
Although there is a potentially large number of variable and constraints,
namely $nK$, we will see later that the degrees of all polynomials involved
here are generally small, and also the number $K$ of sub-intervals is not large.
So, the problem (\ref{SOPE-C_app}) does not have an excessively high complexity.

\begin{remark} \rm
For some polynomial $\varepsilon_i(t) \ge 0$, denote
\be
P_{k,\varepsilon}(t) = \sum_{i=1}^n [p_i(t) - \gamma_{ik}(t) - \varepsilon_i(t)] \hat b_{ik}(t)
+ [\gamma_{ik}(t) + \varepsilon_i(t)] \check b_{ik}(t).
\label{Pnosigneps}
\ee
We note that $P_{k,\varepsilon}(t) \le P_k(t)$.

So, when solving an optimization problem with the constraint $P_k(t) \ge 0$,
the polynomials $\gamma_{ik}(t)$ will generally tend to 
take their smallest possible values that are allowed by the constraints (\ref{picond});
the underlying reason is that the least conservative is the constraint $P_k(t) \ge 0$,
the larger the feasibility domain of the optimization problem and hence
a possibly better solution.
So, if $p_i(t) \le 0$, then it results that $\gamma_{ik}(t)=0$,
and (assuming $\deg \gamma_{ik} \ge \deg p_i$) if $p_i(t) \ge 0$, then $\gamma_{ik}(t)=p_i(t)$;
the approximation (\ref{Pnosign}) actually coincides
(a posteriori, after the optimization problem is solved and $p_i(t)$ is available)
with that for known signs (\ref{Psign}).
If $p_i(t)$ changes the sign on the current interval, then $\gamma_{ik}(t)$ is
the best upper approximation to $\max(p_i(t), 0)$.

The above values of $\gamma_{ik}(t)$ are reached only if the constraint
$P_k(t) \ge 0$ is active on the interval $[t_{k-1}, t_k]$.
However, the important conclusion is that the construction (\ref{Pnosign})-(\ref{picond})
naturally gives the best lower approximation of $f(t)$ by a polynomial $P_k(t)$,
given the bounding polynomials (\ref{approxint}).
\QED
\end{remark}

\section{One-sided polynomial approximations of a set of exponentials}
\label{sec:expapp}

We present here the tools for finding the polynomial approximation (\ref{approxint})
and the intervals (\ref{interv}) for a set of exponentials
$e^{-\lambda_i t}$, $i=1:n$.
We also suggest how this can be done for other functions.

\subsection{Approximation of a single exponential}

We start by showing how to find lower and upper polynomial approximations to
a single exponential $e^{-\lambda t}$ on an interval $[\tau_0, \tau_1]$.
This is made via \cite[Th.4]{BDV66}, taking advantage of the fact that
the derivatives of all orders of the exponential have constant sign.

We present in detail only one case, that of the lower approximation $\check b(t)$
of odd degree $\nu = 2 \ell - 1$, to a function $\phi(t)$ whose derivative of order $\nu+1$
is nonnegative on $[\tau_0, \tau_1]$ (which is the case of our exponential).
The approximation is optimal in the sense that the norm
\be
\int_{\tau_0}^{\tau_1} w(t) [\phi(t) - \check b(t)] dt
\label{anorm}
\ee
is minimized, where $w(t)$ is a weight function.
Denote by $x_1$, \ldots, $x_\ell$ the zeros of the $\ell$-th order
polynomial from the sequence of polynomials which are orthogonal on
$[\tau_0, \tau_1]$ with respect to $w(t)$.
Then, the coefficients of $\check b(t)$ can be found by solving the linear system
given by the equations
\be
\check b(x_k) = \phi(x_k), \ \check b'(x_k) = \phi'(x_k), \ \ k=1:\ell.
\label{eqb}
\ee
So, if the zeros of $w(t)$ are readily available, then the computation
of $\check b(t)$ is very simple and effective.

The other cases, corresponding to the three other combinations of the
parity of $\nu$ and the sign of the $\nu+1$-th derivative of $\phi(t)$
are similar, but involve the zeros of the polynomials that are orthogonal
with respect to $(t-\tau_0) w(t)$, $(\tau_1-t) w(t)$ and $(t-\tau_0)(\tau_1-t) w(t)$,
and one or both of the interval ends.

For the sake of quick computation, we chose the weight
$w(t) = \sqrt{(\tau_1 - t)(t - \tau_0)}$, which generates Chebyshev polynomials
of the first kind.
Their roots are
\be
x_k = \frac{\tau_1 - \tau_0}{2} \cos \frac{(2k -1) \pi}{2 \ell} 
+ \frac{\tau_1 + \tau_0}{2}, \ \ k=1:\ell.
\label{rootscheb}
\ee
(The other three weights above generate Chebyshev polynomials of the second kind
and Jacobi $(-\frac{1}{2}, \frac{1}{2})$ and $(\frac{1}{2}, -\frac{1}{2})$ polynomials,
whose roots are also availabe via simple formulas.)

\begin{figure}
  \centering
  \bt{ccc}
  \includegraphics[scale=0.35]{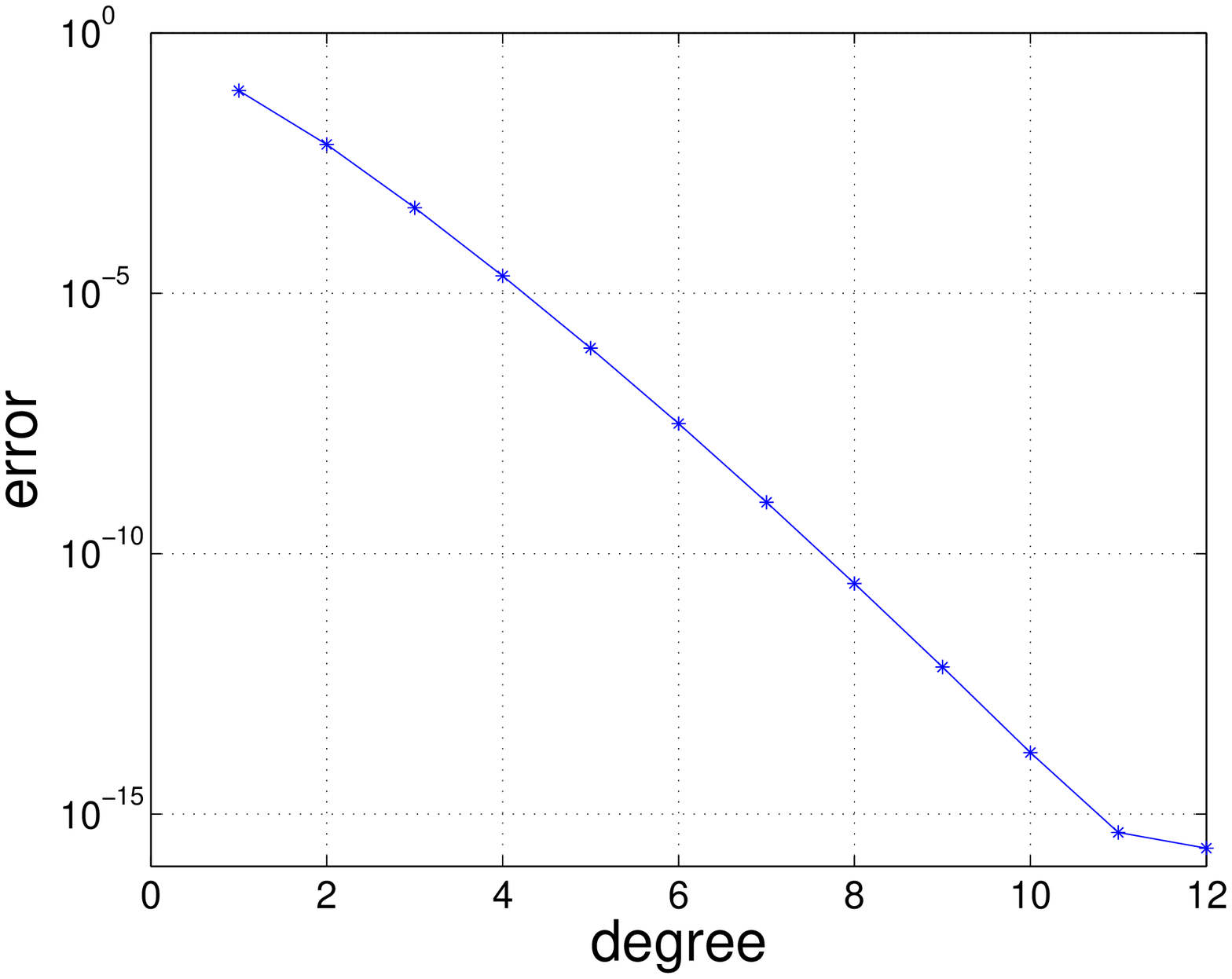}
  & &
  \includegraphics[scale=0.35]{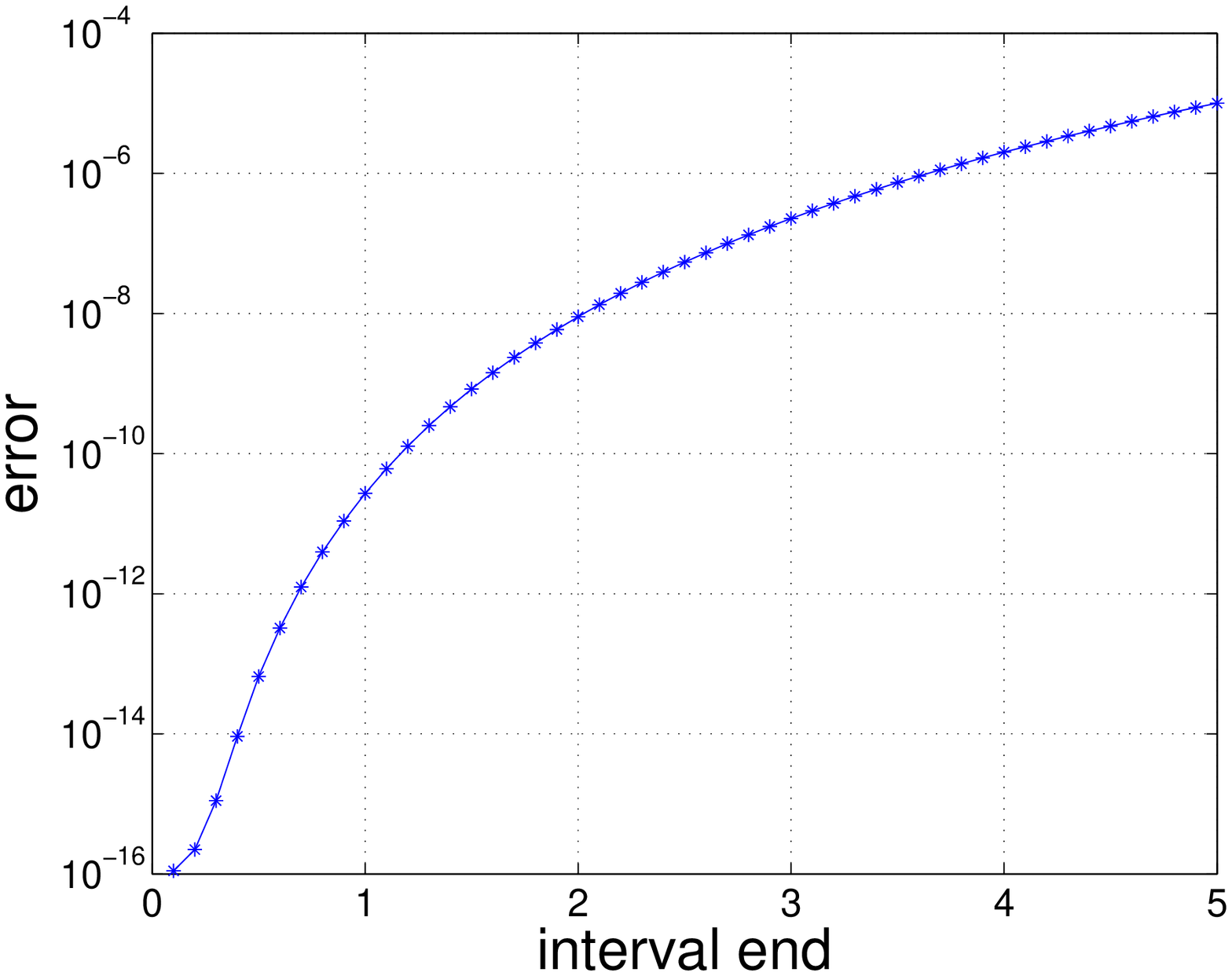}
  \et
  \caption{One-sided polynomial approximation error for $e^{-t}$.
  Left: for various degrees on the interval $[0,1]$.
  Right: for degree 8 and intervals from 0 to the value on the horizontal axis.}
  \label{fig:approx_err}
\end{figure}

\begin{example} \rm
To have an idea of the approximation error, we plot in Figure \ref{fig:approx_err}
the maximum value of
$\hat b(t) - e^{-t}$ and  $e^{-t} - \check b(t)$ for several degrees of
the polynomial and intervals starting from 0 and ending in various points
up to 5.
One can see that, for example, a polynomial of degree 8 gives an error
smaller than $10^{-10}$ for intervals included in $[0, 1.1]$.
The approximation becomes unreliable when the approximation error
approaches $10^{-16}$, in the sense that it may be no longer one-sided.
However, errors of order $10^{-12}$ appear perfectly obtainable.
\QED
\end{example}

Although the infinity norm would be preferrable to an integral norm like (\ref{anorm}),
in order to minimize the maximum approximation error, the computation
of such a norm would involve a Remez-like algorithm
which significantly increases the computation time with relatively small benefits
(at least in the case of the exponential and our choices of intervals).

An alternative to the above construction was employed in \cite{ShWa01} and
works for general functions, not only in the restrictive conditions from \cite{BDV66}.
Using a grid like (\ref{rootscheb}), one builds Chebyshev interpolation polynomials,
for whom an approximation error bound is available.
Subtracting or adding the bound value to the polynomial gives lower or upper
one-sided approximations.
A more refined approach \cite{PaTr09}, based on a sound numerical implementation
of the Remez algorithm, can give nearly optimal solutions.

\subsection{Approximations of a set of exponentials}

We go now to our full problem, approximating a set of exponentials on a reunion
of intervals.
We assume a target approximation error $\epsilon$.
Given the exponentials $e^{-\lambda_i t}$, $i=1:n$, the problem is how
to find the intervals (\ref{interv}) and the degrees $\nu_{ik}$ of the polynomials
(\ref{approxint}) such that
\be
\ba{lcl}
\ba{l}
\max_{t \in [t_{k-1}, t_k]} |\hat b_{ik}(t) - e^{-\lambda_i t}| \le \epsilon, \\
\max_{t \in [t_{k-1}, t_k]} |e^{-\lambda_i t} - \check b_{ik}(t)| \le \epsilon,
\ea
& & \forall i=1:n, \ \forall k = 1:K.
\ea
\label{approxerr}
\ee
The solution is obviously not unique, but we aim to find one that is
computationally advantageous.
The exponential function allows a very cheap solution.

\begin{remark} \rm
\label{rem:expapprox}
Assume that the maximum approximation error of a polynomial
of degree $\nu$ approximating $e^{-t}$ on the interval $[0, \tau]$ is $\varepsilon$.
Then, there exist polynomials of degree $\nu$ such that:

\begin{itemize}
 \item $e^{-\lambda t}$ is approximated on $[0, \tau / \lambda]$ with error
at most $\varepsilon$;
 \item $e^{-t}$ is approximated on $[x, x+\tau]$ with error
at most $\varepsilon e^{-x}$.
\end{itemize}

Combining these two results, it follows that $e^{-\lambda t}$ is approximated on 
$[x, x+\tau/\lambda]$ with error at most $\varepsilon e^{-\lambda x}$.
\QED
\end{remark}

\begin{example} \rm
Let us put the problem differently.
We take $\epsilon = 10^{-10}$ and $\nu = 8$.
We want to approximate $e^{-3t}$ on an interval $[1, 1+\theta]$
as large as possible, with approximation error less than $\epsilon$.
From the above Remark, the error for $e^{-t}$ on $[0, 3 \theta]$
should be less than $e^3 \epsilon \approx 2 \cdot 10^{-9}$.
Looking in Figure \ref{fig:approx_err} (right), we see that
this happens for $3 \theta \approx 1.6$.
Indeed, taking $\theta = 0.55$ gives the desired accuracy.
\QED
\end{example}

The approximation procedure we propose is as follows.
Besides the input data listed in the beginning of this section,
we assume that a maximum degree $\nu_{\max}$ is given.

1. We make a table $T(\tau, \nu)$ of approximation errors given by
polynomials built as described in this section for approximating
$e^{-t}$ on $[0,\tau]$.
The errors are measured on a grid;
there can be a few tens of $\tau$ values
and $\nu$ can go from 0 to 10 or 12.
(The result would be a "cartesian product" of the two graphs
from Figure \ref{fig:approx_err}.)
This table, having a few hundred entries, is built a single time.
Even so, the computation time for a $50 \times 12$ table was only 0.25 seconds
on a standard desktop, which is negligible.

2. Assume that we have found $t_{k-1}$ and we search $t_k$.
For each $\lambda_i$, using the table and Remark \ref{rem:expapprox}, we find
the interval length $d_{ik}$ such that the approximation error for
$e^{- \lambda_i t}$ on $[t_{k-1}, t_{k-1} + d_{ik}]$ is at most $\epsilon$.
More precisely, we seek in the table the value $\tau$ for which
$T(\tau, \nu_{\max})$ has the largest value smaller than
$\epsilon \cdot e^{ \lambda_i t_{k-1}}$ and set
$d_{ik} = \tau / \lambda_i$.
Finally, to ensure that the approximation error is respected for all
exponentials, we take the smallest interval length and put
$t_k = t_{k-1} + \min_{i=1:n} d_{ik}$.
This iterative procedure, ending when $t_k \ge t_f$, is extremely fast since
it involves only table searches.


For the sake of numerical accuracy, we optionally can reduce the degrees
of the approximating polynomials for the exponentials that are not
deciding the length of an interval.
(Note that the highest degree of a polynomial actually decides the complexity.)
Denoting $d_k = t_k - t_{k-1}$, for each $i$ for which $d_{ik} > d_k$
we search in the table the smallest degree $\nu_{ik}$ for which
$T(\tilde \tau, \nu) e^{\lambda_i \tau_{k-1}}$ is smaller than $\epsilon$,
where $\tilde \tau$ is the smallest grid value larger than $d_k/\lambda_i$.
Remark \ref{rem:expapprox} ensures that the error made by approximating
$e^{- \lambda_i t}$ on $[t_{k-1}, t_k]$ with a polynomial of degree $\nu_{ik}$
is at most $\epsilon$.

\begin{example} \rm
Let us approximate the exponentials with exponents $0.3$, 1 and 3
on the interval $[0,10]$, with error less then $\epsilon = 10^{-10}$,
using polynomials of degree at most $\nu_{\max} = 8$.
The table is built with step value $0.2$ for $\tau$.
Using the above procedure, we find that 10 sub-intervals are necessary.
The first is $[0,0.33]$, on which the degrees of the approximating
polynomials are 5, 6 and 8; the fastest decaying exponential
naturally needs the highest degree.
The last sub-interval is $[9.27,10]$ and the degrees are 6, 5 and 1.
They are smaller than $\nu_{\max}$ because the sub-interval is cut
short by $t_f=10$; for example, the previous sub-interval is $[6.67,9.27]$;
since the exponentials are decaying, the sub-intervals are longer as
$t$ grows.
Note that now the slowest exponential sets the degree;
the fastest has almost vanished.
\QED
\end{example}

Also for improving numerical accuracy, it is useful to move the
approximation on an interval centered in the origin.
Similarly to a Vandermonde system, the linear system (\ref{eqb})
tends to become ill-conditioned when the roots (\ref{rootscheb})
have all the same sign and large absolute value.
So, denoting $\theta_k = t_{k-1} + d_k/2$,
instead of working with $p_i(t) e^{-\lambda_i t}$ for $t \in [t_{k-1}, t_k]$,
we work with $p_i(t + \theta_k) e^{-\lambda_i \theta_k} e^{-\lambda_i t}$
for $t \in [-d_k/2, d_k/2]$.
Since usually the degree $\degpi_i$ of $p_i(t)$ is small,
the bad numerical effect of computing the coefficients of $p_i(t + \theta_k)$
(as a polynomial in $t$)
is much smaller than the reduction of the condition number of the
system solved for approximating $e^{-\lambda_i t}$ on $[-d_k/2, d_k/2]$.

It is clear that the above procedure is particularly fit for the exponential function.
For other functions, the polynomial approximations must be computed explicitly
for each sub-interval, in a trial-and-error procedure for finding the right
polynomial degrees and the sub-intervals lengths; this is essentially
what Chebfun does, but in our case more flexibility is allowed since we
don't (and cannot) aim to approximation within machine precision.

It is also obvious that imposing an approximation error $\epsilon$
for each exponential, does not make the approximation error in (\ref{fpos})
be of the same size.
However, an a posteriori analysis of the function $f(t)$, in particular of the
values of the computed polynomials $p_i(t)$, can help estimate the
actual error.
A single new run with a smaller $\epsilon$ would be required, since
the new optimized function should be near the previous one
and hence the worst-case accuracy becomes predictable.

\section{Fitting an empirical density}
\label{sec:fit}

Let us now come back to our prototype problem (\ref{Jsope}) of
fitting an empirical density and consider first the simpler SOE case (\ref{soe}).
We start by studying a helpful particular case.


\subsection{Exponents in arithmetic progression}
\label{sec:aritprog}

Let us assume that the exponents are known and form an arithmetic progression,
meaning that
\be
\lambda_i = \lambda_1 + (i-1) q, \ \ i=1:n,
\label{aritprog}
\ee
where $q > 0$ is the ratio.
In this case, the function (\ref{soe}) has the form
\be
f(t) = e^{-\lambda_1 t} \sum_{i=1}^n  \alpha_i e^{-(i-1)qt}.
\label{soepol}
\ee
Denoting $x = e^{-qt}$, the condition (\ref{fpos}) becomes
\be
\sum_{i=1}^n  \alpha_i x^{i-1} \ge 0, \ \ \forall x \in [e^{-q t_f}, e^{-q t_0}].
\label{fposx}
\ee
This is a polynomial positivity condition, easy to impose through LMIs \cite{Nest00}.
Hence the SOPE-C problem (\ref{Jsope}) is equivalent to an SDP problem.
No approximation is required, the problem is naturally polynomial. 

The method from \cite{Duf07} works with models based on the progression (\ref{aritprog})
and transforms $f(t)$ into a polynomial as above.
However, positivity is imposed by quite rudimentary means
(fitting a SOE to the square root of the target pdf, then squaring, which
gives also a SOE, but with more terms).
Moreover, the values $\lambda_1$ and $q$ defining the arithmetic progression
are found by just trying different values until a satisfying result is obtained.  

\subsection{An algorithm for the SOE-G problem}

The arithmetic progression case can be used for initialization in an
iterative procedure for the general problem SOE-G.
The iterative part uses the polynomial approximation
described in Sections \ref{sec:appol} and \ref{sec:expapp}.

{1. \em Initialization using sparse arithmetic progression.}
Let us assume that an estimate  $\tilde \lambda_1$ of $\lambda_1$ is available.
This can be obtained like in \cite{SOH12}, using the tail of $h(t)$;
for large values of $t$, the slowest exponential dominates the others.
Alternatively, we can just take $\tilde \lambda_1$ sufficiently small,
since it is enough to have $\tilde \lambda_1 < \lambda_1$.

We then attempt to solve SOE-G by making a more general assumption than
in Section \ref{sec:aritprog} and searching exponents estimations $\tilde \lambda_i$
that {\em belong} to an arithmetic progression
\be
\tilde \lambda_i = \tilde \lambda_1 + \mu_i q,
\label{lamap}
\ee
where $\mu_i$ are unknown positive integers.

To this purpose, we work with the function
\be
\tilde f(t) = \sum_{i=0}^N  \tilde \alpha_i e^{-(\tilde \lambda_1 + qi) t}
= e^{- \tilde \lambda_1 t} \sum_{i=0}^N  \tilde \alpha_i e^{-qit},
\label{fp}
\ee
with given $q$ and $N$.
In principle, $q$ should be small enough to cover decently well
the possible intervals where the exponents lie and $N$ should be
as large as computationally acceptable
(e.g.\ 100 is certainly good, but one can consider going to 200 and even beyond).
Since we seek a sparse solution, with only $n$ nonzero coefficients $\tilde \alpha_i$,
we modify the problem (\ref{Jsope}) by adding a sparsity-promoting term to the
criterion and transforming it into
\be
\ba{cl}
 \min & J(\tilde f) + \beta \sum_{i=0}^N \varpi_i |\tilde \alpha_i| \\
 \mbox{s.t.} & \tilde f(t) \ge 0, \ \forall t \in [0,\infty)
\ea
\label{fpopt}
\ee
where $\beta$ and $\varpi_i$, $i=0:N$, are weighting constants.
The second term of the criterion is a weighted 1-norm of the coefficients
vector $\bm{\tilde \alpha}_i$, denoted $\| \bm{\tilde \alpha} \|_{\varpi,1}$;
the most meaningful choice appears to be $\varpi_i = 1/\tilde \lambda_i$,
which takes into account that $\int_0^\infty e^{-\lambda t} = 1/\lambda$
and thus implicitly normalizes the exponentials from (\ref{soe}).

The problem (\ref{fpopt}) is convex in the coefficients $\tilde \alpha_i$.
The positivity constraint can be expressed as the positivity of a polynomial $\tilde f(x)$
on $[0,1]$, by substituting $x = e^{-qt}$ as in Section \ref{sec:aritprog}.
If the constant $\beta$ is large enough, many of the coefficients
$\tilde \alpha_i$ will be small and only few will have significant values.
We take the largest $n$ of them (in weighted absolute value $\varpi_i |\tilde \alpha_i|$)
and the corresponding $\lambda_i$ are given by their positions.

We can then re-solve (\ref{fpopt}) by imposing that only the chosen
$n$ coefficients are nonzero, finding thus the optimal coefficients $\alpha_i$
for the selected $\lambda_i$
(this can have the advantage of exact positivity constraint via polynomials,
but can be skipped by going directly to the iterative step).

The problem (\ref{fpopt}) could be replaced with a minimization of
$\| \bm{\tilde \alpha} \|_{\varpi,1}$
with a bounded $J(\tilde f)$; in this case, the bound should be chosen instead of $\beta$;
this may make more sense if some value of the criterion is already available.

\smallskip
{2. \em Iterative part.}
With the above initialization, we can start the iterative part of the
algorithm, where the coefficients and the exponents are optimized alternatively.

For given exponents $\lambda_i$, we optimize the coefficients $\alpha_i$ by solving
the approximation (\ref{SOPE-C_app}) of the SOE-C problem (\ref{Jsope}),
as shown in Sections \ref{sec:appol} and \ref{sec:expapp}.
The auxiliary variables $\gamma_{ik}(t)$ are just scalars.

For optimizing the exponents $\lambda_i$ we attempt small gradient steps.
The gradient is
\be
\frac{\partial J(f)}{\partial \lambda_i} = - \frac{2 \alpha_i}{M} 
\sum_{m=1}^M w_m t_m [ f(t_m) - h_m ] e^{- \lambda_i t_m}.
\label{Jsoegrad}
\ee
The step size is $\varsigma$, so the new values of the exponents are
\be
\lambda_i \leftarrow \lambda_i - \varsigma \frac{\partial J(f)}{\partial \lambda_i}.
\label{gradstep}
\ee
(If the resulting $\lambda_i$ would become negative, we reduce the
step size such that they stay positive.)
Then, the SOE-C problem is (approximately) solved with the new exponents $\lambda_i$.
If the value of the criterion does not decrease, we restore the previous
exponents, halve the step size $\varsigma$ and recompute (\ref{gradstep}).
Note that by modifying $\lambda_i$ we usually improve the criterion but
may go out of the positivity domain.
Solving (\ref{SOPE-C_app}) means returning back to it, but not
necessarily with a better criterion value.

The iterations continue as long as the improvement is significant
or the step size is not very small.
Of course, there is no guarantee of optimality, but our main purpose
is to illustrate the kernel of our approach---the polynomial approximation idea.
Even with such a simple descent procedure, the results are satisfactory
in the test problems we will report in Section \ref{sec:res}.

\subsection{Extension to SOPE}

The iterative part of the above procedure can be used for true SOPE functions
(\ref{sope}), with some $\degpi_i > 0$, with the minor change of the
gradient expression (\ref{Jsoegrad}) into
\be
\frac{\partial J(f)}{\partial \lambda_i} = - \frac{2}{M} 
\sum_{m=1}^M w_m t_m p_i(t_m) [ f(t_m) - h_m ] e^{- \lambda_i t_m}.
\label{Jsopegrad}
\ee

However, the arithmetic progression trick used in the SOE case is
no more possible, since the transformation $x=e^{-qt}$ no longer
leads to a polynomial.
Instead, we simply use the exponents resulting from a SOE solution
(not necessarily fully optimized) as initialization for the SOPE problem.
As confirmed by numerical evidence, this seems effective especially
when the polynomials degrees $\degpi_i$ from (\ref{sope}) are small.

\section{Results}
\label{sec:res}

The polynomial approximation method described in this paper
has been implemented using POS3POLY \cite{SiDu13oe} for CVX \cite{CVX}
and can be downloaded from \url{http://www.schur.pub.ro/sope},
together with the programs solving the problems presented in this section.
The simplest description of the constraint (\ref{fpos}), with a SOE function
(\ref{soe}) characterized by the variable vector of coefficients {\tt alpha}
and the constant vector of exponents {\tt lambda} is 
\begin{verbatim}
   alpha == pos_soe( lambda, t0, tf );
\end{verbatim}
where {\tt t0} and {\tt tf} are the positivity interval ends.
This is consistent with CVX style and the parameters of the approximation
are hidden, although the user can control them if so desired.

We start illustrating the behavior of our method by solving a problem
proposed in \cite{SOH12}, that is slightly different in nature from (\ref{Jsope}),
but finally has the same form.

\begin{figure}
  \centering
  \includegraphics[scale=0.4]{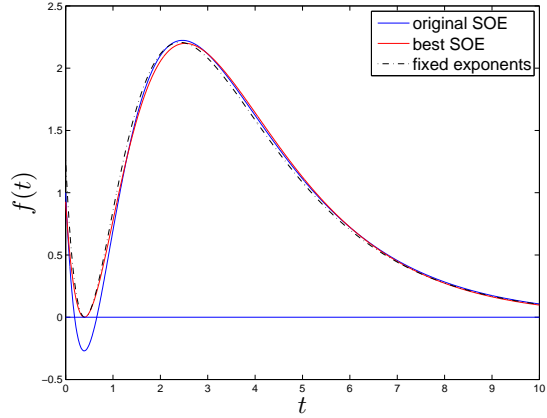}\\
  \caption{Graphs of the functions from Example \ref{ex:sexton}.
   Blue: target $h(t)$ (\ref{sexton}).
   Red: our best approximation (\ref{bestsext}).
   Black, dashed: best approximation (\ref{bestfixexp}),
   with the same exponents as $h(t)$.}
  \label{fig:sexton}
\end{figure}

\begin{example} \rm
\label{ex:sexton}
Unlike an empirical pdf, the SOE
\be
h(t)=16e^{-\frac{t}{2}}-30e^{-t}+15e^{-2t}
\label{sexton}
\ee
has also negative values.
We want to find the nearest positive SOE (\ref{soe}) to $h(t)$,
by minimizing
\be
J(f) = \int_0^\infty [f(t) - h(t)]^2 dt.
\label{Jsexton}
\ee
This criterion is convex quadratic in $\alpha_i$, like (\ref{Jsope}),
but is rational in $\lambda_i$.
The expression of the gradient with respect to $\lambda_i$ is omitted here,
for brevity, but easy to obtain.

Using the tail of $h(t)$ like in \cite{SOH12}, the smallest exponent
estimation is $\tilde \lambda_1 = 0.4873$.
With $q=0.01$, $N=100$ and $\beta = 0.0005$ in (\ref{fp}--\ref{fpopt}),
the initialization step of our algorithm gives 
$\lambda_1 = 0.5073$, $\lambda_2 = 1.0773$, $\lambda_3 = 1.4873$
and a criterion value $J = 0.0433$.
Positivity is imposed on the interval $[0,10]$, which is large enough
to ensure it on $[0, \infty]$.

The iterative step improves it to $J=0.0425$, the result being
\be
f(t)=19.91 e^{- 0.5315 t} - 48.63 e^{- 1.0264 t} + 29.66 e^{- 1.5177 t}.
\label{bestsext}
\ee
Figure \ref{fig:sexton} presents the graph of this function and of the target $h(t)$.
The whole design, containing 30 iterations, took less than 1 minute on
a standard desktop computer.
The polynomial approximations of the exponentials were made with
$\nu_{\max} = 8$ and $\epsilon = 10^{-8}$.
The minimum value of the SOE (\ref{bestsext}), computed on
a very fine grid, is $3.5 \cdot 10^{-7}$.
Taking a smaller $\epsilon$ reduces this value only if the
accuracy of the SDP solver is also increased.
For example, setting {\tt cvx\_precision best} (which means that CVX
iterates as long as the criterion can be decreased without numerical trouble),
a value of $\epsilon = 10^{-12}$ (which is still safe for the polynomial
approximation) leads to a minimum of the SOE of $3.6 \cdot 10^{-11}$.

In the above case, the initialization step already gave a good result.
However, the iterative step can significantly decrease the criterion
if the approximation is poor.
For example, with $q=0.03$, the initialization step gives $J=0.2227$
and the final result is $J=0.0439$.

Solving the SOE-C problem with the exponents from (\ref{sexton}) kept fixed,
the optimal criterion is 0.0712.
The same value is given by a dedicated algorithm that takes into account
that the exponents are actually part of an arithmetic progression and
imposes positivity via polynomials; the resulting optimal SOE, whose graph
is also shown in Figure \ref{fig:sexton}, is
\be
f(t) = 15.5243 e^{-\frac{t}{2}} - 28.5073 e^{-t} + 14.2410 e^{-2t}.
\label{bestfixexp}
\ee
The result reported in \cite{SOH12} is
\be
f(t)=16e^{-\frac{t}{2}}-29.946e^{-t}+ 15.5385e^{-2t}.
\ee
Somewhat surprisingly, the corresponding criterion is $J=0.0933$.
\QED
\end{example}

\medskip
Many methods in the literature are tested using standard densities,
and our next two examples will be of this type.
Given a continuous pdf $h(t)$, we simply discretize it on $M=100$ equidistant
values in the interval $[0,t_f]$ and optimize the criterion (\ref{Jsope})
through the approximation (\ref{SOPE-C_app}).
The values $\nu_{\max} = 8$ and $\epsilon = 10^{-8}$ are kept throughout
the rest of the paper.
The other parameters may change values, but we will not report them,
since they can be found in the Matlab files.

\begin{figure}
  \centering
  \includegraphics[scale=0.4]{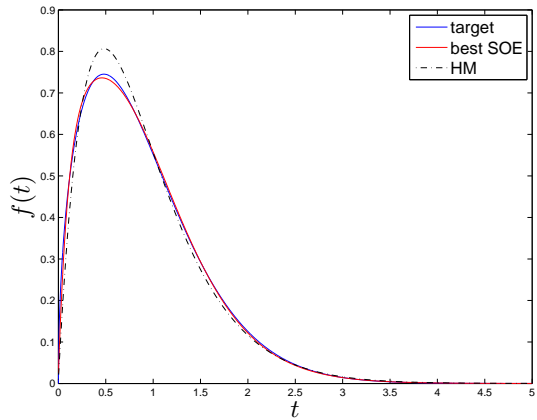}\\
  \caption{Graphs of the functions from Example \ref{ex:hm}.
   Blue: target $h(t)$ (\ref{W1}).
   Red: our best approximation (\ref{W1sol}).
   Black, dashed: the approximation from \cite{HaMa98}.}
  \label{fig:weibull}
\end{figure}

\begin{example} \rm
\label{ex:hm}
We discuss in detail the example W1 from \cite{HaMa98}, where the target
is the Weibull pdf
\be
h(t) = \frac{k}{\eta} \left( \frac{t}{\eta} \right)^{k-1} e^{(-t/\eta)^k},
\label{W1}
\ee
with $\eta=1$, $k=1.5$.
We approximate it with a SOE with $n=4$ terms on the interval $[0,5]$.
The solution
\be
f(t) = 130.12 e^{-2.7535 t} - 329.39 e^{-3.1707 t} + 228.96 e^{-3.5850 t}
- 29.665 e^{-4.7543 t}
\label{W1sol}
\ee
gives $J = 4.82 \cdot 10^{-5}$.
The result is illustrated in Figure \ref{fig:weibull}, together
with the best 4-term SOE reported in \cite{HaMa98}
(where the approximation is made on the Laplace transform),
for which the criterion is only $J=5.90 \cdot 10^{-4}$.
Also, our result is visually similar to the order six model
from \cite[Fig.\ 2]{ANO96}.

For other examples from \cite{HaMa98} involving the Weibull and lognormal
distribution, where we could compare the solutions only graphically
(since the result was given only in this form),
we have obtained better results with the same number of exponentials
or similar results with fewer exponentials.
Note also that the method from \cite{HaMa98} allows complex values for the exponents;
however, for example W1 the exponents were real.
\QED
\end{example}

\begin{example} \rm
\label{ex:dufre}
We turn now to some examples from \cite{Duf07}, where the 
complementary cumulative density function (ccdf)
\be
H(t) = \int_t^\infty h(\tau) d \tau
\label{ccdf}
\ee
was used instead of the pdf $h(t)$ and modeled with a SOE (\ref{soe}).
However, since the derivative of a SOE $F(t)$ is also a SOE $f(t)$,
we have solved a problem very similar to (\ref{Jsope}), namely
\be
\ba{ccl} \min & & \ds J(F) = \frac{1}{M} \sum_{m=1}^M w_m [ F(t_m) - H_m ]^2 \\
\mbox{s.t.} & & f(t) \ge 0, \ \forall t \in [t_0, t_f] \\
            & & F'(t) = - f(t) \\ 
            & & F(t_0) = 1, \ F(t_f) \ge 0
\ea
\label{JsopeF}
\ee
The last two constraints, that do not appear in (\ref{Jsope}), are
linear in the coefficients of the SOE and hence
the algorithms described in this paper can be immediately adapted.

The two ccdf discussed in this example are those of the Pareto
\be
H(t) = \frac{1}{1+t}
\ee
and lognormal
\be
H(t) = 1 - \frac{1}{2} \text{erfc} \left(-\frac{\log{(t-\mu)}}{\sigma\sqrt{2}} \right)
\ee
distributions.
For the lognormal, the parameters are $\mu = 0$, $\sigma = 0.5$.
The intervals on which optimization is performed are $[0,20]$ for Pareto
and $[0,6]$ for the lognormal.

We present in Table \ref{tab:dufre} the values of the criteria (\ref{JsopeF})
for several values of the number $n$ of exponentials and three methods.
The first is an optimized version of the method from \cite{Duf07}, which
used exponents in arithmetic progression.
The first exponent and the progression ratio were found by trial and error;
we use the same values as reported there.
In \cite{Duf07}, the SOE is approximated using Jacobi polynomials, but
with no explicit optimization.
We use instead the exact optimization method sketched in section \ref{sec:aritprog}.
Hence we obtain better results than those reported in \cite{Duf07}.
(Note that the maximum error is reported in \cite{Duf07}, while
we optimize an LS criterion and still get better maximum error.)

The other two methods are ours.
In addition to SOE, we report now also some results for SOPE models.
The values shown in Table \ref{tab:dufre} are obtained with
degrees $\degpi_i = 1$, $i=1:n$, in (\ref{sope}).
For Pareto, when comparing the number $\degpi = n+\sum_{i=1}^n \degpi_i$
of coefficients $\alpha_{i,j}$,
the SOPE model is sometimes more efficient than SOE.
For example, the criterion for SOPE with $n=5$ and hence $\degpi=10$
is smaller than that for SOE with $n=\degpi=11$.
Also, if we take $\degpi=5$ via $\degpi_1 = 0$, $\degpi_2 = \degpi_3 = 1$,
the criterion is $J = 5.53 \cdot 10^{-6}$, smaller than the value for SOE and $n=5$.
For the lognormal (as well as for the previous examples, where SOPE was
not mentioned), SOE appears more adequate than SOPE.

Comparing now our SOE method with the least-squares optimal version
of the method from \cite{Duf07}, we notice two distinct behaviors.
For small $n$, our method is significantly better, showing that
the search for good exponents is successful and that using an arithmetic
progression may be quite restrictive.
However, for larger $n$, the methods give comparable results.
Now the arithmetic progression seems a satisfactory model.
Numerical accuracy properties of the SDP solver may also come into play,
since the fit is actually very good and the values of the criterion very small.
Finally, we note that our method, by its nature, targets small values of $n$,
where the parameters of the model may have a significance.
Taking exponents in a long arithmetic progression resembles more
what was called non-parametric method in the introduction.
\QED
\end{example}

\begin{table}
  \centering
  \begin{tabular}{|c|c|c|c|c|c|}
  \hline
Target & Approximation &  \multicolumn{4}{|c|}{$n$} \\
 \cline{3-6}
ccdf & & 3 & 5 & 11 & 21 \\
  \hline
        & \cite{Duf07}, LS optimal & $1.84 \cdot 10^{-2}$ & $5.02 \cdot 10^{-4}$ & $1.40 \cdot 10^{-6}$ & $1.52 \cdot 10^{-8}$ \\
Pareto  & Our SOE             & $1.28 \cdot 10^{-4}$ & $2.36 \cdot 10^{-5}$ & $1.73 \cdot 10^{-6}$ & $3.42 \cdot 10^{-9}$ \\ 
        & Our SOPE            & $2.87 \cdot 10^{-6}$ & $1.27 \cdot 10^{-7}$ & $1.33 \cdot 10^{-6}$ & $4.73 \cdot 10^{-9}$ \\
\hline
           & \cite{Duf07}, LS optimal & $8.97 \cdot 10^{-3}$ & $1.86 \cdot 10^{-3}$ & $3.01 \cdot 10^{-6}$ & $4.50 \cdot 10^{-7}$ \\
lognormal  & Our SOE             & $1.39 \cdot 10^{-3}$ & $6.68 \cdot 10^{-5}$ & $3.83 \cdot 10^{-6}$ & $3.75 \cdot 10^{-6}$ \\
           & Our SOPE            & $1.29 \cdot 10^{-4}$ & $7.14 \cdot 10^{-6}$ & $4.00 \cdot 10^{-6}$ & $5.26 \cdot 10^{-6}$ \\
\hline
\end{tabular}
  \caption{Values of the criterion (\ref{JsopeF}) for Example \ref{ex:dufre}.}
  \label{tab:dufre}
\end{table}

\begin{figure}
  \centering
  \includegraphics[scale=0.4]{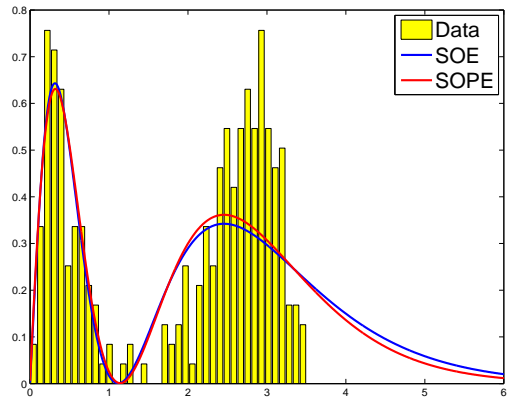}\\
  \caption{Graphs of the functions from Example \ref{ex:oldfaith}.}
  \label{fig:oldfaith}
\end{figure}

\begin{example} \rm
\label{ex:oldfaith}
We fit now SOE and SOPE models to data belonging to a popular time series,
the eruption durations of the Old Faithful geyser.
The problem is notoriously difficult since the pdf appears to
have two disjoint intervals as support.
There are 272 values in the series, which are grouped into $M=40$ equally spaced bins
and shifted towards the origin with 1.6 minutes, which is the smallest duration
of an eruption.
Figure \ref{fig:oldfaith} shows the resulting histogram and two approximations
with $\degpi=6$ parameters, a SOE with $n=6$ and a SOPE with $n=4$ and
$\degpi_1 = 0$, $\degpi_2 = \degpi_3 = 1$, $\degpi_4 = 0$.
The optimization is performed on the interval $[0,10]$;
zero values are appended to the histogram for times greater
than the longer duration.
The values of the criterion are $0.0112$ for SOE and $0.0104$ for SOPE,
so again the SOPE model is more appropriate.
The first peak of the data and the valley are quite well followed,
while the approximation is worse for the second peak.

A visual comparison with Fig.\ 4b from \cite{KimT11} shows that our models
with only 6 parameters have a better fit than some (older) methods with
order 10 models and are only slightly worse than the model proposed in \cite{KimT11}.
All the models investigated there are more general than SOPE.
Anyway, our purpose in this example is not to prove that SOPE is a good
model for the Old Faithful data, but that our method is flexible enough
to give reasonably good results in this case.
Better approximation of the data can be obtained by shortening the optimization
interval, but then, like in \cite{KimT11}, the pdf is not strictly decreasing
after the second peak, but has a third small peak.
\QED
\end{example}

\section{Conclusion and future work}

We have presented a method for solving SIP problems with inequality constraints
involving functions (\ref{sope}) that are a sums of exponentials multiplied with
polynomials with variable coefficients.
The method is based on one-sided polynomial approximations of the exponentials
on sub-intervals, which allow transformation of SIP into SDP and hence
reliable computation of near-optimal solutions that are guaranteed to
respect the constraints.
We showed through several examples that our method is computationally
attractive and gives good results compared to methods based on different
principles.

Future work will be dedicated to the investigation of approximations for
other functions. The exponential has properties that allowed us several
tricks easing the computation of the polynomial approximation.
Extension to general functions is not trivial, but is certainly possible.
We aim to build an user-friendly library for convex SIP based on the
principles set up in this paper.


\end{document}